\documentstyle[12pt]{article}
\topmargin by -\headsep

\evensidemargin \oddsidemargin \marginparwidth 0.5in

\begin{document}

\topmargin -.6in

\def\rf#1{(\ref{eq:#1})}
\def\lab#1{\label{eq:#1}}
\def\nonu{\nonumber}
\def\br{\begin{eqnarray}}
\def\er{\end{eqnarray}}
\def\be{\begin{equation}}
\def\ee{\end{equation}}
\def\eq{\!\!\!\! &=& \!\!\!\! }
\def\foot#1{\footnotemark\footnotetext{#1}}
\def\lb{\lbrack}
\def\rb{\rbrack}
\def\llangle{\left\langle}
\def\rrangle{\right\rangle}
\def\blangle{\Bigl\langle}
\def\brangle{\Bigr\rangle}
\def\llb{\left\lbrack}
\def\rrb{\right\rbrack}
\def\Blb{\Bigl\lbrack}
\def\Brb{\Bigr\rbrack}
\def\lcurl{\left\{}
\def\rcurl{\right\}}
\def\({\left(}
\def\){\right)}
\def\v{\vert}                     %% vertical bars
\def\bv{\bigm\vert}               %%
\def\Bgv{\;\Bigg\vert}            %%
\def\bgv{\bigg\vert}              %%
\def\lskip{\vskip\baselineskip\vskip-\parskip\noindent}
\def\mskp{\par\vskip 0.3cm \par\noindent}
\def\sskp{\par\vskip 0.15cm \par\noindent}
\relax

\def\tr{\mathop{\rm tr}}                  % tr - small trace
\def\Tr{\mathop{\rm Tr}}                  % Tr - big trace
\newcommand\partder[2]{{{\partial {#1}}\over{\partial {#2}}}}
                                                  % partial derivative
\newcommand\funcder[2]{{{\delta {#1}}\over{\delta {#2}}}}
                                                % functional derivative
\newcommand\Bil[2]{\Bigl\langle {#1} \Bigg\vert {#2} \Bigr\rangle}  %% <.|.>
\newcommand\bil[2]{\left\langle {#1} \bigg\vert {#2} \right\rangle} %% <.|.>
\newcommand\me[2]{\langle {#1}\vert {#2} \rangle} %% <.|.>

\newcommand\sbr[2]{\left\lbrack\,{#1}\, ,\,{#2}\,\right\rbrack} % commutator
\newcommand\Sbr[2]{\Bigl\lbrack\,{#1}\, ,\,{#2}\,\Bigr\rbrack} % commutator
%%(Large)
\newcommand\pbr[2]{\{\,{#1}\, ,\,{#2}\,\}}       % Poisson brackets
\newcommand\Pbr[2]{\Bigl\{ \,{#1}\, ,\,{#2}\,\Bigr\}}  % Poisson brackets
%%(large)
\newcommand\pbbr[2]{\lcurl\,{#1}\, ,\,{#2}\,\rcurl}  % Poisson brackets
%%(left-right)
\newcommand\sumi[1]{\sum_{#1}^{\infty}}   %% summation till infinity

\def\a{\alpha}
\def\b{\beta}
\def\c{\chi}
\def\d{\delta}
\def\D{\Delta}
\def\eps{\epsilon}
\def\vareps{\varepsilon}
\def\g{\gamma}
\def\G{\Gamma}
\def\grad{\nabla}
\def\h{{1\over 2}}
\def\k{\kappa}
\def\l{\lambda}
\def\L{\Lambda}
\def\m{\mu}
\def\n{\nu}
\def\om{\omega}
\def\O{\Omega}
\def\p{\phi}
\def\P{\Phi}
\def\pa{\partial}
\def\pr{\prime}
\def\ra{\rightarrow}
\def\lra{\longrightarrow}
\def\s{\sigma}
\def\S{\Sigma}
\def\t{\tau}
\def\th{\theta}
\def\Th{\Theta}
\def\z{\zeta}
\def\ti{\tilde}
\def\wti{\widetilde}
\def\one{\hbox{{1}\kern-.25em\hbox{l}}}

\def\cA{{\cal A}}
\def\cB{{\cal B}}
\def\cC{{\cal C}}
\def\cD{{\cal D}}
\def\cE{{\cal E}}
\def\cH{{\cal H}}
\def\cJ{{\cal J}}
\def\cL{{\cal L}}
\def\cM{{\cal M}}
\def\cN{{\cal N}}
\def\cP{{\cal P}}
\def\cQ{{\cal Q}}
\def\cR{{\cal R}}
\def\cS{{\cal S}}
\def\cT{{\cal T}}
\def\cU{{\cal U}}
\def\cV{{\cal V}}
\def\cW{{\cal W}}
\def\cY{{\cal Y}}

\def\symp#1{{\cal S}{\cal D}if\!\! f \, ({#1})}
\def\esymp#1{{\wti {\cal S}{\cal D}if\!\! f} \, ({#1})}
\def\Symp#1{{\rm SDiff}\, ({#1})}
\def\eSymp#1{{\wti {\rm SDiff}}\, ({#1})}
\def\vol#1{{{\cal D}if\!\! f}_0 ({#1})}
\def\Vol#1{{\rm Diff}_0 ({#1})}

\def\Winf{{\bf W_\infty}}               % Linear W-infinity
\def\Win1{{\bf W_{1+\infty}}}           % Linear W-1+infinity
\def\nWinf{{\bf {\hat W}_\infty}}       % Nonlinear W-infinity
\def\PsDA{\Psi{\cal DO}}
                   % algebra of all pseudo-differential operators
%%
\newcommand{\nit}{\noindent}
\newcommand{\ct}[1]{\cite{#1}}
\newcommand{\bi}[1]{\bibitem{#1}}

\begin{titlepage}
\vspace*{-1cm} \phantom{bla}

\begin{center}{\large \textbf{Infinite Deformed Groups and Their \\
Geometrical and Physical Applications}} \foot{This article
represent the report presented at the \emph{Bogolyubov Kyiv
conference "Modern Problems of Mathematics and Theoretical
Physics"} held at Kyiv (Ukraine), 13-17 September 2004}
\end{center} \vskip .3in

\begin{center}
\textsc{Serhiy E. Samokhvalov} \vskip 2mm {\small
\emph{State Technical University, Dniprodzerzhinsk, Ukraine}\\
\centerline{\emph{e-mail: samokhval@dstu.dp.ua}}}
\end{center}
\vskip .3in

\begin{abstract}

In this article we give the realization of the Klein's Program for
geometrical structures (Riemannian spaces and fiber bundles with
connection) with arbitrary variable curvature within the framework
of infinite deformed groups. These groups generalize gauge groups
to the case of nontrivial action on the base space of bundles with
use of idea of groups deformations.

We also show that infinite deformed groups give a group-theoretic
description of gauge fields (gravitational field with its metric
or vierbein part similarly to gauge fields of internal symmetry)
which is alternative for their geometrical interpretation.

\end{abstract}

\end{titlepage}

\centerline{\textbf{Introduction}} \vskip 3mm{\small

An infinite (local gauge) symmetry lays in a basis of modern
theories of fundamental interactions. Theory of gravitation
(general relativity) is based on the idea of covariance with
respect to the group of space-time diffeomorphisms. Theories of
strong and electroweak interactions are gauge theories of internal
symmetry. Moreover, the existence of these interactions is
considered to be necessary for ensuring the local gauge
symmetries.

But any physical theory can be written in covariant form without
introduction of a gravitational field. Similarly, as first had
been emphasized in [1], for any theory with global internal
symmetry $G$ corresponding gauge symmetry $G^{g}$ can be ensured
without introduction of nontrivial gauge fields with help of pure
gauging. Presence of the gravitational or the gauge fields of
internal symmetry is being manifested in presence of deformation -
curvature of Riemannian space, or fiber bundles with connection
accordingly.

Formally nontrivial gauge fields are being entered by continuation
of derivatives up to covariant derivatives $\partial_{\mu}
\rightarrow \nabla_{\mu}$. Their commutators characterize strength
of a field, which is considered, and from the geometrical point of
view - curvature of corresponding space. On the other hand,
covariant derivatives set infinitesimal space-time translations in
gauge fields (curved spaces). That is why it is possible to
suppose, that for introduction of nontrivial arbitrary gauge
fields it is necessary to consider groups (wider than gauge groups
$G^{g}$ ) which would generalize gauge groups $G^{g}$ to a case of
\emph{nontrivial action on space-time manifold and contain the
information about arbitrary gauge fields in which occur a motion}.
Consequently, such groups must \emph{contain the information about
appropriate geometrical structures with arbitrary variable
curvature} and set these geometrical structures on manifolds where
they act. Hence from the mathematical point of view such groups
should realize the \emph{Klein's Erlanger Program} [2] for these
geometrical structures.

For a long time it was considered that such groups do not exist.
\emph{E.Cartan} [3] has named the situation in question as
\emph{Riemann-Klein's antagonism} - antagonism between
\emph{Riemann's} and \emph{Klein's} approach to geometry. There
are attempts of modifying of the Klein's Program for geometrical
structures with arbitrary variable curvature by means of refusal
of group structure of used transformations with usage of
categories [4], quasigroups [5] and so on. One can encounter
widespread opinion that \emph{nonassociativity} is an algebraic
equivalent of the geometrical notion of a curvature [5].

In our article we shall show that \emph{realization of the Klein's
Program for the geometrical structures with arbitrary variable
curvature (Riemannian space and fiber bundles with connection) can
be executed within the framework of the so called infinite
deformed groups which generalize of gauge groups to the case of
nontrivial action on the base space of bundles with use of idea of
groups deformations}.

Such groups have been constructed in [6]. Klein's Erlanger Program
was realized for fiber bundles with connection in [7] and for
Riemannian space in [8].

This is important for physics because the widely known gauge
approaches to gravity (see, for example, [9]) in fact gives gauge
interpretation neither to metric fields nor to vierbein ones. An
interpretation of these as connections in appropriate fibrings has
been achieved in a way of introduction (explicitly or implicitly)
of the assumption about existence of the background flat space
(see, for example, [10]). That is unnatural for gravity. The
reason for these difficulties lies in the fact that the fiber
bundles formalism is appropriate only for the internal symmetry
Lie group, which do not act on the space-time manifold. But for
the gravity this restriction is obviously meaningless because it
does not permit consider gravity as the gauge theory of the
translation group.

In this article we also shall show that \emph{infinite deformed
groups give a group-theoretic description of gauge fields
(gravitational field with its metric or vierbein part similarly to
gauge fields of internal symmetry) which is alternative for their
geometrical interpretation} [11].

This approach allow to overcome the well known Coleman-Mandula
no-go theorem within the framework of infinite deformed groups and
gives new possibilities to unification gravity with gauge theories
of internal symmetry [12].

\vskip 10mm

\centerline{\textbf{1. Generalized Gauge Groups}} \vskip 3mm
{\small

Gauge groups of internal symmetry $G^{g}$ are a special case of
infinite groups and have simple group structure - the infinite
direct product of the finite-parameter Lie groups
$G^{g}=\prod_{x\in M} G$ where product takes on all points $x$ of
the space-time manifold $M$. Groups $G^{g}$ act on $M$ trivially:
$x'^{\mu}=x^{\mu}$.

For the aim of a generalization of groups $G^{g}$ to the case of
nontrivial action on the space-time manifold $M$, let's now
consider a Lie group $G_{M}$ with coordinates $\tilde{g}^{\alpha}$
(indexes $\alpha, \beta, \gamma, \delta$) and the multiplication
law: \vskip 5mm

\centerline{$(\tilde{g}\cdot\tilde{g}')^{\alpha}=\tilde{\varphi}^{\alpha}
(\tilde{g},\tilde{g}')$,} \vskip 3mm

\noindent which act (perhaps inefficiency) on the space-time
manifold $M$ with coordinates $x^{\mu}$ (indexes $\mu, \nu, \pi,
\rho, \sigma$) according to the formula: \vskip 5mm

\centerline{$x'^{\mu}=\tilde{f}^{\mu}(x, \tilde{g})$.}\vskip 3mm

The infinite Lie group $G^{g}_{M}$ is parameterized by smooth
functions $\tilde{g}^{\alpha}(x)$ which meet the condition: \vskip
5mm

\centerline{$\det\{d_{\nu}\tilde{f}^{\mu}(x,
\tilde{g}(x))\}\neq0$, $\forall x\in M$,}\vskip 3mm

\noindent where $d_{\nu}:=d/dx^{\nu}$. The multiplication law in
$G^{g}_{M}$ is determined by the formulae\,:

\begin{equation}
\label{1}
(\tilde{g}\times\tilde{g}')^{\alpha}=\tilde{\varphi}^{\alpha}(\tilde{g}(x),
\tilde{g}'(x')),\,
\end{equation}

\begin{equation}
\label{2} x'^{\mu}=\tilde{f}^{\mu}(x, \tilde{g}(x)).
\end{equation}

\noindent It is a simple matter to check that these operations
truly make $G^{g}_{M}$ a group. Formula (2) sets the action of
$G^{g}_{M}$ on $M$. In the case of trivial action of the group
$G_{M}$ on $M$, $G^{g}_{M}$ becomes the ordinary gauge group
$G^{g}=\prod_{x\in M} G$. The groups $G^{g}_{M}$ we call
\emph{generalized gauge groups}.

For the clearing of the groups deformations idea we shall consider
spheres of different radius $R$. All of them have isomorphic
isometry groups - groups of rotations $O(3)$. The information
about radius of the spheres is in structural constants of groups
$O(3)$, which in the certain coordinates may be written as:
$F^3_{12}=1 / R^2$, $F^2_{13}=-1$, $F^1_{23}=1$. Isomorphisms of
groups $O(3)$, which change $R$, correspond to deformations of the
sphere.

For gauge groups $G^g_M$ some isomorphisms also play a role of
deformations of space of groups representations, but as against
isomorphisms of finite-parameter Lie groups such isomorphisms are
more substantial, as these allow to independently deform space in
its different points.

Let us pass from the group $G^g_M$ to the group $G^{gH}_M$
isomorphic to it by the formula $g^a (x)=H^a (x, \tilde{g}(x))$
(Latin indexes assume the same values as the corresponding Greek
ones). The smooth functions $H^a (x, \tilde{g})$ have the
properties:\vskip 5mm

 1)\ $H^{a}(x,0)=0$\quad$\forall x \in M$;\vskip 3mm

 2)\ $\exists K^a (x,g)$:\quad$K^a (x,H(x,\tilde g))=
 \tilde g^\alpha$\quad$\forall\tilde g \in G,\ x \in M$.\vskip 3mm

\noindent The group $G^{gH}_M$ multiplication law is determined by
its isomorphism to the group $G^{g}_M$ and formulae (1), (2):

\begin{equation}
\label{3} (g \ast g')^a (x)=\varphi^a (x,g(x),g'(x')) := H^a
(x,\tilde\varphi(K(x,g(x)),K(x',g'(x')))),
\end{equation}

\begin{equation}
\label{4} x'^{\mu}=f^{\mu}(x,g(x)) := \tilde f^\mu (x,K(x,g(x))).
\end{equation}

\noindent Formula (4) sets the action of $G^{gH}_M$ on $M$.

Such transformations between the groups $G^{g}_M$ and $G^{gH}_M$
we call \emph{deformations} of infinite Lie groups, since
(together with changing of the multiplication law) the
corresponding deformations of geometric structures of manifolds
subjected to group action are directly associated with them. The
functions $H^a (x, \tilde g)$ we call \emph{deformation
functions}, functions \vskip 5mm

\centerline{$h(x)^a_\alpha := \partial H^a (x,\tilde g) /\partial
\tilde g^\alpha |_{\tilde g = 0}$}\vskip 3mm

\noindent - \emph{deformation coefficients}, and the groups
$G^{gH}_M$ - \emph{infinite (generalized gauge) deformed groups}.

Let us consider expansion:

\begin{equation}
\label{5} \varphi^a (x,g,g')=g^a + g'^a + \gamma(x)^a_{bc}\ g^b
g'^c + \frac{1}{2} {\rho(x)^a}_{bcd}\ g^d g'^b g'^c + \ldots
\end{equation}

\noindent The functions $\varphi^a$, setting the multiplication
law (3) in the group $G^{gH}_M$, are explicitly $x$-dependent, so
the coefficients of expansion (5) are $x$-dependent as well. So,
$x$-dependent became structure coefficients of group $G^{gH}_M$
(\emph{structure functions} versus structure constants for
ordinary Lie groups)

\begin{equation}
\label{6} F(x)^a_{bc} := \gamma(x)^a_{bc} - \gamma(x)^a_{cd}
\end{equation}

\noindent and coefficients

\begin{equation}
\label{7} R(x)^a_{dbc} := \rho(x)^a_{dbc} - \rho(x)^a_{dcb} ,
\end{equation}

\noindent which we call \emph{curvature coefficients} of the
deformed group $G^{gH}_M$.

Since \vskip 5mm

\centerline{$\xi(x)_a^\mu := \partial_a f^\mu_H (x,g) |_{g = 0} =
h(x)^\alpha_a \tilde \xi (x)^\mu_\alpha$, }\vskip 3mm

\noindent where $\partial_b :=
\partial / \partial g^b$ and $h(x)^\alpha_a$ is reciprocal to the
$h(x)_\alpha^a$ matrix, the generators $X_a := \xi(x)_a^\mu
\partial_\mu$ ($\partial_\mu := \partial / \partial x^\mu$) of the
deformed group $G^{gH}_M$ are expressed through the generators
$\tilde X_\alpha := \tilde \xi(x)_\alpha^\mu \partial_\mu$ of the
group $G^{g}_M$ with the help of deformation coefficients: \vskip
5mm

\centerline{$X_a = h(x)_a^\alpha \tilde X_\alpha$. }\vskip 3mm

\noindent So, in infinitesimal (algebraic) level, deformation is
reduced to independent in every point $x \in M$ nondegenerate
liner transformations of generators of the initial Lie group.
\vskip 3mm

\textbf{Theorem 1.} \emph{Commutators of generators of the
infinite (generalized gauge) deformed group are liner combinations
of generators with structure functions as coefficients} [6]:

\begin{equation}
\label{8} [X_a , X_b ] = F(x)^c_{ab} X_c .
\end{equation} \vskip 3mm

For the generalized gauge nondeformed group $G^{g}_M$ we have:
\vskip 5mm

\centerline{$[\tilde X_\alpha , \tilde X_\beta ] = \tilde
F^\gamma_{\alpha\beta} \tilde X_\gamma$, }\vskip 3mm

\noindent where $\tilde F^\gamma_{\alpha\beta}$ - structure
constants of the initial Lie group $G_M$.

\vskip 10mm}

\centerline{\textbf{2. Group-Theoretic Description of Connections
in Fiber Bundles}} \centerline{\textbf{and Gauge Fields of
Internal Symmetry}} \vskip 3mm {\small

Let $P=M \times V$ be a principal bundle with the base $M$
(space-time) and a structure group $V$ with coordinates $\tilde
\upsilon^i$ (indexes $i, j, k$) and the multiplication law
$(\tilde \upsilon \cdot \tilde \upsilon')^i = \tilde \varphi^i
(\tilde \upsilon, \tilde \upsilon')$. As usually, we define the
left $l_{\tilde \upsilon} : P = M \times V \rightarrow P' = M
\times \tilde \upsilon^{-1} \cdot V$ and the right $r_{\tilde
\upsilon} : P = M \times V \rightarrow P' = M \times V \cdot
\tilde \upsilon$ action $V$ on $P$.

Let's consider a group $G_M = T_M \otimes V$ where $T_M$ is the
group of space-time translations. The group $G_M$ is parameterized
by the pair $\tilde t^\mu$, $\tilde \upsilon^i$, has the
multiplication law: \vskip 5mm

\centerline{$(\tilde g \cdot \tilde g')^\mu = \tilde t^\mu +
\tilde t'^\mu,\quad(\tilde g \cdot \tilde g')^i = \tilde \varphi^i
(\tilde \upsilon, \tilde \upsilon')$ }\vskip 3mm

\noindent and act on the $M$ inefficiently: $x'^\mu = x^\mu +
\tilde t^\mu$. One can define the left action of the group $G_M$
on the principal bundle $P$ : \vskip 5mm

\centerline{$x'^\mu = x^\mu + \tilde t^\mu,\quad\upsilon'^i =
l^i_{\tilde \upsilon} (\upsilon)$. }\vskip 3mm

The group $G^g_M$ is parameterized by the functions $\tilde
t^\mu(x)$, $\tilde \upsilon^i (x)$ which meet the condition:
\vskip 5mm

\centerline{$\det\{\delta_{\nu}^{\mu}+\partial_{\nu}\tilde{t}^{\mu}(x)\}\neq0$,
$\forall x\in M$.}\vskip 3mm

\noindent The multiplication law in $G^g_M$ is:

\begin{equation}
\label{9} (\tilde g \times \tilde g')^\mu (x) = \tilde t^\mu (x) +
\tilde t'^\mu (x'),\quad(\tilde g \times \tilde g')^i (x) = \tilde
\varphi^i (\tilde \upsilon (x), \tilde \upsilon' (x')) ,
\end{equation}

\begin{equation}
\label{10} x'^\mu = x^\mu + \tilde t^\mu (x) ,
\end{equation}

\noindent where (10) determines the inefficient action of $G^g_M$
on $M$ with the kern of inefficiency - gauge group $V^g$. The
group $G^g_M$ has the structure $Diff\ M \times ) V^g$, act on $P$
as:\vskip 5mm

\centerline{$x'^\mu = x^\mu + \tilde t^\mu (x),\quad\upsilon'^i =
l^i_{\tilde \upsilon (x)} (\upsilon) $}\vskip 3mm

\noindent and is the group $aut\ P$ of  automorphisms of the
principal bundle $P$.

Let us deform the group $G^g_M \rightarrow G^{gH}_M$ with help of
deformation functions with additional properties:\vskip 5mm

3)\ $H^{\mu}(x, \tilde t, \tilde \upsilon) = \tilde
t^\mu$\quad$\forall \tilde t \in T,\ \tilde \upsilon \in V,\ x \in
M$;\vskip 3mm

4)\ $H^{i}(x, 0, \tilde \upsilon) = \tilde
\upsilon^i$\quad$\forall \tilde \upsilon \in V,\ x \in M$.\vskip
3mm

\noindent The deformed group $G^{gH}_M$ is parameterized by the
functions:\vskip 5mm

\centerline{$t^\mu (x) = \tilde t^\mu (x),\quad \upsilon^i (x) =
H^i (x,\tilde t (x), \tilde \upsilon (x)).$}\vskip 3mm

\noindent Obviously, the group $G^{gH}_M$, as well as the group
$G^{g}_M$, has the structure $Diff\ M \times ) V^g$ and act on $P$
as:

\begin{equation}
\label{11} x'^\mu = x^\mu + t^\mu (x) ,\quad\upsilon'^i =
l^i_{K(x,t(x),\upsilon(x))} (\upsilon) ,
\end{equation}

\noindent where functions $K^i(x,t(x),\upsilon(x))$ are determined
by equation: $K^i(x,t(x),\upsilon(x)) = \tilde \upsilon^i (x)$.
Properties 3), 4) result in the fact that among deformation
coefficients of the group $G^{gH}_M$, $x$-dependent is only
$h(x)^i_\mu = \partial_{\tilde\mu} H^i (x, \tilde t, 0) |_{\tilde
t =0} =: -A(x)^i_\mu$ (where $\partial_{\tilde \mu} := \partial /
\partial \tilde t^\mu)$.

Generators of the $G^{gH}_M$-action on $P$ (11) are split in the
pair: \vskip 5mm

\centerline{$X_{\mu}=\partial_{\mu}+A(x)_{\mu}^{i}\tilde{X}_{i},\quad
X_{i}=\tilde{X}_{i}$} \vskip 3mm

\noindent where $\tilde{X}_{i}$ are generators of the left action
of the group $V$ on $P$. This results in the natural splitting of
tangent spaces $T_{u}$ in any point $u \in P$ into the direct sum
$T_{u}=T\tau_{u}\oplus T\upsilon_{u}$ subspaces: \vskip 5mm

\centerline{$T\tau_{u}=\{t^{\mu}X_{\mu}\},\quad
T\upsilon_{u}=\{\upsilon^{i}X_{i}\}$.}\vskip 3mm

\noindent The distribution $T\tau_{u}$ is invariant with respect
to the right action of the group $V$ on $P$, and $T\upsilon_{u}$
is tangent to the fiber. So $T\tau_{u}$ one can treated as
horizontal subspaces of $T_{u}$ and generators $X_{\mu}$ - as
covariant derivatives. This set in the principal bundle $P$
connection and deformation coefficients $A(x)_{\mu}^{i}$ are the
coordinates of the connection form, which on submanifold $M
\subset P$ may be written as $\omega^{i}=A(x)_{\mu}^{i}dx^{\mu}$.
Necessary condition of existence of group $G_{M}^{gH}$ (8) for
generators $X_{\mu}$ yield:

\begin{equation}
\label{12} [X_{\mu},X_{\nu}]=F(x)_{\mu\nu}^{i}X_{i} ,
\end{equation}

\noindent where

\begin{equation}
\label{13}
F(x)_{\mu\nu}^{i}X_{i}=\tilde{F}_{jk}^{i}A(x)_{\mu}^{j}A(x)_{\nu}^{k}
+\partial_{\mu}A(x)_{\nu}^{i}-\partial_{\nu}A(x)_{\mu}^{i}
\end{equation}

\noindent - structure functions of the group $G_{M}^{gH}$ and
$\tilde{F}_{jk}^{i}$ - structure constants of the Lie group $V$.
Relationship (12) one can write as:

\begin{equation}
\label{14}
d\omega^{i}=-\frac{1}{2}\tilde{F}_{jk}^{i}\omega^{j}\wedge
\omega^{k}+\Omega^{i},
\end{equation}

\noindent where form \vskip 5mm

\centerline{$\Omega^{i}=\frac{1}{2}
F(x)_{\mu\nu}^{i}dx^{\mu}\wedge dx^{\nu}$} \vskip 3mm

\noindent play the role of the curvature form on submanifold $M$.
So, equation (14) is a structure equation for connection, which
has be set on the principal bundle $P$ by action of group
$G_{M}^{gH}$. \vskip 3mm

\textbf{Theorem 2.} \emph{Acting on the principal bundle $P=M
\times V$ deformed group $G_{M}^{gH}=Diff\ M \times ) V^g$ sets on
$P$ structure of connection. Any connection on the principal
bundle $P=M \times V$ may be set thus} [7]. \vskip 3mm

This theorem realizes Klein's Erlanger Program for fiber bundles
$P=M \times V$ with connection.

We should emphasize that for the setting of a geometrical
structure on $P$ it is enough to consider the infinitesimal action
(11) of the group $G_{M}^{gH}$.

The potentials of gauge fields of internal symmetry are identified
with deformation coefficients $A(x)_{\mu}^{i}$, a strength tensor
- with structure functions $F(x)_{\mu\nu}^{i}$ of the group
$G_{M}^{gH}$. All groups $G_{M}^{gH}$ obtained one from another by
internal automorphisms, which are generated by the elements
$\upsilon(x)\in V^{g}$, describe the same gauge field. These
automorphisms lead to gauge transformations for fields
$A(x)_{\mu}^{i}$ and for infinitesimal $\upsilon^{i}(x)$ yield:

\begin{equation}
\label{15}
A'(x)_{\mu}^{i}=A(x)_{\mu}^{i}-\tilde{F}_{jk}^{i}A(x)_{\mu}^{j}\upsilon^{k}(x)
-\partial_{\mu}\upsilon^{i}(x).
\end{equation}

\vskip 10mm

\centerline{\textbf{3. Group-Theoretic Description of Riemannian
Spaces}} \centerline{\textbf{and Gravitational Fields}} \vskip 3mm
{\small

The structure of Riemannian space is a special case of structure
of affine connection in vierbein bundle and consequently it can be
set by the way described above with the application of the
deformed group $Diff M\times)SO(n)^{g}$. If we force the
coordination of connection with a metric and vanishing of torsion,
generators of translations
$X_{\mu}=\partial_{\mu}+\Gamma(x)_{\mu}^{(mn)}\tilde{S}_{(mn)}$,
where $\tilde{S}_{(mn)}$ is generators of group $SO(n)$, become
covariant derivatives in Riemannian space. For the setting of
Riemannian structure by such means, it is enough to consider the
group $Diff M\times)SO(n)^{g}$ on algebraic level - on level its
generators.

Potentials of a gravitational field in the given approach are
represented by the connection coefficients
$\Gamma(x)_{\mu}^{(mn)}$ instead of the metrics or verbein fields
that would correspond to sense of a gravitational field as a gauge
field of translation group which is born by an energy - momentum
tensor, instead of spin.

Now we shall show that the Riemannian structure on $M$ is
naturally set also by a narrower group than $Diff
M\times)SO(n)^{g}$, namely, the deformed group of diffeomorphisms
$T_{M}^{gH}$, though it demands consideration of its action on $M$
up to the second order on parameters.

Let $G_{M}=T_{M}$ where $T_{M}$ is the group of space-time
translations. In this case
$(\tilde{t}\cdot\tilde{t'})^{\mu}=\tilde{t}^{\mu}+\tilde{t'}^{\mu}$
and $x'^{\mu}=x^{\mu}+ \tilde t^\mu$. The group is parameterized
by the functions $\tilde{t}^{\mu}(x)$, which meet the condition:
\vskip 5mm

\centerline{$\det\{\delta_{\nu}^{\mu}+\partial_{\nu}\tilde{t}^{\mu}(x)\}\neq0$,
$\forall x\in M$.}\vskip 3mm

\noindent The multiplication law in $T_{M}^{g}$ is:

\begin{equation}
\label{16} (\tilde t \times \tilde t')^\mu (x) = \tilde t^\mu (x)
+ \tilde t'^\mu (x'),
\end{equation}

\begin{equation}
\label{17} x'^\mu = x^\mu + \tilde t^\mu (x) ,
\end{equation}

\noindent where (17) determines the action of $T_{M}^{g}$ on $M$.
The multiplication law indicates that $T_{M}^{g}$ is the group of
space-time diffeomorphisms $DiffM$ in additive parameterization.
Thus, in the approach considered, the group $T_{M}^{g}=DiffM$
becomes the gauge group of local translations. The generators of
the $T_{M}^{g}$-action on $M$ (17) are simply derivatives
$\tilde{X}_{\mu}=\partial_{\mu}$ and this fact corresponds to the
case of the flat space $M$ and the absence of gravitational field.

Let us deform the group $T_{M}^{g}\rightarrow T_{M}^{gH}$:
$t^m(x)=H^m(x,\tilde{t}(x))$. The multiplication law in
$T_{M}^{gH}$ is determined by the formulae:

\begin{equation}
\label{18} (t \ast t')^m (x)=\varphi^m (x,t(x),t'(x')) :=
H^m(x,K(x,t(x)+K(x',t'(x'))),
\end{equation}

\begin{equation}
\label{19} x'^{\mu}=f^{\mu}(x,t(x)) := x^\mu + K^{\mu}(x,t(x)).
\end{equation}

\noindent Formula (19) sets the action of $T_{M}^{gH}$ on $M$.

Let us consider expansion:

\begin{equation}
\label{20} H^{m}(x,\tilde{t})=h(x)_{\mu}^{m}[\tilde{t}^{\mu}+
\frac{1}{2}\Gamma(x)_{\nu\rho}^{\mu}\tilde{t}^{\nu}\tilde{t}^{\rho}+
\frac{1}{6}\Delta(x)_{\nu\rho\sigma}^{\mu}\tilde{t}^{\nu}\tilde{t}^{\rho}
\tilde{t}^{\sigma}].
\end{equation}

\noindent With usage of the formula (18), for coefficients of
expansion (5) we can obtain:

\begin{equation}
\label{21} \gamma_{pn}^{m}=h_{\mu}^{m}(\Gamma_{pn}^{\mu}+
h_{p}^{\nu}\partial_{\nu}h_{n}^{\mu}),
\end{equation}

\begin{equation}
\label{22} \rho_{prn}^{m}=h_{\mu}^{m}(\Delta_{prn}^{\mu}-
\Gamma_{ns}^{\mu}\Gamma_{pr}^{s}-
h_{n}^{\nu}\partial_{\nu}\Gamma_{\pi\rho}^{\mu}
h_{p}^{\pi}h_{r}^{\rho}).
\end{equation}

\noindent So, formulae (6), (7) for structure coefficients and
curvature coefficients of deformed group $T_{M}^{gH}$ yield:

\begin{equation}
\label{23}
F_{\mu\nu}^{n}=-(\partial_{\mu}h_{\nu}^{n}-\partial_{\nu}h_{\mu}^{n}),
\end{equation}

\begin{equation}
\label{24} R_{\rho\pi\nu}^{\mu}=
\partial_{\pi}\Gamma_{\nu\rho}^{\mu}-\partial_{\nu}\Gamma_{\pi\rho}^{\mu}
+\Gamma_{\pi\sigma}^{\mu}\Gamma_{\nu\rho}^{\sigma}-
\Gamma_{\nu\sigma}^{\mu}\Gamma_{\pi\rho}^{\sigma}.
\end{equation}

\noindent In this formulae deformation coefficients
$h(x)_{\mu}^{m}$ and $h(x)_{m}^{\mu}$ we use for changing Greek
index to Latin (and vice versa).

Formulae (23) and (24) say that groups $T_{M}^{gH}$ contain
information about geometrical structure of space $M$ where they
act. The generators $X_{n}=h_{n}^{\nu}\partial_{\nu}$ of the
$T_{M}^{gH}$-action (19) on $M$ can be treated as vierbeins.
Structure functions $F_{\mu\nu}^{n}$ differ from non-holonomity
coefficients only by a factor -1/2.

Let us write the multiplication law of the group $T_{M}^{gH}$ (18)
for infinitesimal second factor:

\begin{equation}
\label{25} (t \ast \tau)^m (x)=t^m
(x)+\lambda(x,t(x))_{n}^{m}\tau^{n}(x'),
\end{equation}

\noindent where
$\lambda(x,t(x))_{n}^{m}:=\partial_{n'}\varphi^{m}(x,t,t')\mid_{t'=0}$.
Formula (25) gives the rule for the addition of vectors, which set
in different points $x$ and $x'$ or a \emph{rule of the parallel
transport} of a vector field $\tau$ from point $x'$ to point
$x$:\vskip 5mm

\centerline{$\tau_{\parallel}^{m}(x)
=\lambda(x,t(x))_{n}^{m}\tau^{n}(x')$,}\vskip 3mm

\noindent or in coordinate basis: \vskip 5mm

\centerline{$\tau_{\parallel}^{\mu}(x)
=\partial_{\tilde{\nu}}H^{\mu}(x,\tilde{t})
\tau^{\nu}(x+\tilde{t})$.}\vskip 3mm

\noindent This formula determines the \emph{covariant derivative}:

\begin{equation}
\label{26} \nabla_{\nu}\tau^{\mu}(x)=\partial_{\nu}\tau^{\mu}(x)+
\Gamma(x)_{\sigma\nu}^{\mu}\tau^{\sigma}(x),
\end{equation}

\noindent where functions $\Gamma(x)_{\sigma\nu}^{\mu}$ set the
second order of expansion (20) and play the role of coefficients
of an affine connection. They are symmetric on the bottom indexes,
so torsion equals zero. If $\eta_{mn}$ is a metric of a flat
space, in the manifold $M$ we can determine metrics
$g_{\mu\nu}=h_{\mu}^{m}h_{\nu}^{n}\eta_{mn}$. If we force
$\gamma_{msn}^{\cdot}+\gamma_{nsm}^{\cdot}=0$ (lowering indexes we
fulfill with help of metric $\eta_{mn}$), we can show that
coefficients $\Gamma_{\mu\nu}^{\rho}$ of expansion (20) may be
written as:

\begin{equation}
\label{27}
\Gamma_{\mu\nu}^{\rho}=\frac{1}{2}g^{\rho\sigma}(\partial_{\mu}g_{\nu\sigma}
+\partial_{\nu}g_{\mu\sigma}-\partial_{\sigma}g_{\mu\nu}).
\end{equation}

\noindent So, these coefficients coincide with \emph{Christoffel
symbols} $\{_{\mu\nu}^{\rho}\}$ and curvature coefficients
$R_{\lambda\kappa\nu}^{\mu}$ of the group $T_{M}^{gH}$ coincide
with the \emph{Riemann curvature tensor}. \vskip 3mm

\textbf{Theorem 3.} \emph{Acting on the manifold $M$ the deformed
group $T_{M}^{gH}$ sets on $M$ structure of a Riemannian space.
Any Riemannian structure on the manifold $M$ may be set thus}
[8].\vskip 3mm

This theorem realizes Klein's Erlanger Program for Riemannian
space $M$.

Information about Christoffel symbols is contained in the
\emph{second order} of expansion (20) of deformation functions and
about curvature - in functions $\rho_{prn}^{m}$, which determine
the \emph{third order} of expansion (5) of the multiplication law
of the group $T_{M}^{gH}$. So, in this approach we need consider
not only infinitesimal - algebraic - level in the group
$T_{M}^{gH}$ (as in previous section), but higher levels, too.

The gravitational field potentials are identified with deformation
coefficients $h(x)_{\mu}^{m}$, a strength tensor of the
gravitational field - with structure functions $F(x)_{mn}^{p}$ of
the group $T_{M}^{gH}$. All groups $T_{M}^{gH}$ obtained one from
another by internal automorphisms describe the same gravitational
field. These automorphisms, which can always be connected with the
coordinate transformations on $M$, lead to a general covariance
transformation law for fields $h(x)_{\mu}^{m}$ and for
infinitesimal $t^{m}(x)$ yield:

\begin{equation}
\label{28}
h'(x)_{\mu}^{m}=h(x)_{\mu}^{m}-F(x)_{np}^{m}h(x)_{\mu}^{n}t^{p}(x)
-\partial_{\mu}t^{m}(x).
\end{equation}

The transformation law (28) is similar to the transformation law
(15) for potentials of gauge fields of internal symmetry and the
only difference consists in the replacement of structure constants
of the finite-parameter Lie group $V$ by structure functions of
the infinite deformed group $T_{M}^{gH}$. This fact permits us to
interpret the group $T_{M}^{gH}$ as the gauge translation group
and the fields $h(x)_{\mu}^{m}$ as the gauge fields of the
translation group.

So, we show that infinite deformed (generalized gauge) groups:

1. \emph{set on manifolds where they act geometrical structures of
fiber bundles with connection or Riemannian spaces with arbitrary
variable curvature};

2. \emph{give a group-theoretic description of gauge fields of
internal symmetry as well as gravitational fields}.

\vskip 10mm

\centerline{\textbf{References}} \vskip 3mm
{\small

\noindent 1. Ogievetski V.I., Polubarinov I.V. \emph{On the
Meaning of the Gauge Invariance}, Nuovo Cimento \textbf{23}
(1962), 173-180.

\noindent 2. Klein F. \emph{Vergleichende Betrachtungen Uber
Neuere Geometrische Forschungen (Erlanger Programm)}, Erlangen,
1872.

\noindent 3. Cartan E. \emph{Group Theory and Geometry}, 1925.

\noindent 4. Sulanke R., Wintgen P. \emph{Differentialgeometrie
und Faserbundel, Veb Deutscher Verlag der Wissenschaften}, Berlin,
1972.

\noindent 5. Sabinin L.V. \emph{Methods of Nonassociativ Algebras
in Differential Geometry}, 1981.

\noindent 6. Samokhvalov S.E. \emph{Group-Theoretic Description of
Gauge Fields}, Theor. Math. Phys. \textbf{76} (1988), 709-720.

\noindent 7. Samokhvalov S.E. \emph{About the Setting of
Connections in Fiber Bundles by the Acting of Infinite Lie
Groups}, Ukrainian Math. J. \textbf{43} (1991), 1599-1603.

\noindent 8. Samokhvalov S.E. \emph{Group-Theoretic Description of
Riemannian Spaces}, Ukrainian Math. J. \textbf{55} (2003),
1238-1248.

\noindent 9. Hehl F.W., Heyde P., Kerlich G.D., Nester J.M.
\emph{General relativity with spin and torsion: fondations and
prospects}, Rev. Mod. Phys. \textbf{48} (1976), 393-416.

\noindent 10. Cho Y.M. \emph{Einstein lagrangian as the
translational Yang-Mills lagrangian}, Phys. Rev. \textbf{14D}
(1976), 2521-2525.

\noindent 11. Konopleva N.P., Popov V.N. \emph{Gauge Fields},
Atomizdat, Moscow, 1980.

\noindent 12. Samokhvalov S.E. \emph{About the Symmetry of the
Electrogravitational Unification}, Problems of Nuclear Phys. and
Cosmic Rays \textbf{35} (1991), 50-58.

}

%\centerline{\textbf{References}} \vskip 3mm

\end{document}